\documentclass[amstex,12pt,russian,amssymb]{article}

\usepackage{mathtext}
\usepackage[cp1251]{inputenc}
\usepackage[T2A]{fontenc}
\usepackage[russian]{babel}
\usepackage[dvips]{graphicx}
\usepackage{amsmath}
\usepackage{amssymb}
\usepackage{amsxtra}
\usepackage{latexsym}
\usepackage{ifthen}

\begin{document}

\newcounter{lemma}
\newcommand{\lemma}{\par \refstepcounter{lemma}%
{\bf Лемма \arabic{lemma}.}}

\newcounter{corollary}
\newcommand{\corollary}{\par \refstepcounter{corollary}%
{\bf Следствие \arabic{corollary}.}}

\newcounter{remark}
\newcommand{\remark}{\par \refstepcounter{remark}%
{\bf Замечание \arabic{remark}.}}

\newcounter{theorem}
\newcommand{\theorem}{\par \refstepcounter{theorem}%
{\bf Теорема \arabic{theorem}.}}

\newcounter{proposition}
\newcommand{\proposition}{\par \refstepcounter{proposition}%
{\bf Предложение \arabic{proposition}.}}

\renewcommand{\refname}{\centerline{\bf Список литературы}}

\newcommand{\proof}{{\it Доказательство.\,\,}}

\noindent УДК 517.5

{\bf Е.А.~Севостьянов, С.А.~Скворцов} (Житомирский государственный
университет имени Ивана Франко)

{\bf Є.О.~Севостьянов, С.О.~Скворцов} (Житомирський державний
університет імені Івана Фран\-ка)

{\bf E.A.~Sevost'yanov, S.A.~Skvortsov} (Zhytomyr Ivan Franko State
University)

\medskip
{\bf О сходимости отображений в метрических пространствах с прямыми
и обратными модульными условиями}

{\bf Про збіжність відображень у метричних просторах з прямим і
оберненими модульними умовами}

{\bf On convergence of mappings in metric spaces with direct and
inverse modulus conditions}

\medskip\medskip
Для отображений метрических пространств, удовлетворяющих одной
оценке модуля семейств кривых, получен результат о нульмерности
предельного отображении. Доказано, что равномерный предел
последовательности таких отображений нульмерен, как только
мажоранта, отвечающая за искажение семейств кривых, имеет конечное
среднее колебание в каждой точке. Кроме того, для одного класса
гомеоморфизмов метрических пространств получены теоремы о
равностепенной непрерывности семейства обратных отображений.

\medskip\medskip
Для відображень метричних просторів, що задовольняють одну оцінку
модуля сімей кривих, отримано результат  про нульвимірність
граничного відображення. Доведено, що рівномірною границею
послідовності вказаних відображень є нульвимірне відображення, як
тільки мажоранта, що відповідає за спотворення сімей кривих, має
скінченне середнє коливання в кожній точці. Крім того, для одного
класу гомеоморфізмів метричних просторів отримано теореми про
одностайну неперервність сімей обернених відображень.

\medskip\medskip
For mappings in metric spaces satisfying one inequality with respect
to modulus of families of curves, there is proved a lightness of the
uniform limit of mappings mentioned above. It is proved that, the
uniform limit of these mappings is light mapping, whenever a
function which corresponds to distortion of families of curves, is
of finite mean oscillation at every point. Besides that, for one
class of homeomorphisms of metric spaces, there are obtained
theorems about equicontinuity of inverse mappings.

\newpage
{\bf 1. Введение.} Работа посвящена изучению ёмкостно-модульной
техники и квазиконформных отображений в метрических пространствах.
Как известно, указанные вопросы активно исследуются последнее время
(см., напр., \cite{AS}, \cite{HK}, \cite{He}, \cite{OR} и
\cite{RSa}). Основная цель статьи -- изучить сходимость отображений
в метрических пространствах, в которых принципиально возможно
применение аппарата модулей семейств кривых. В статье поочерёдно
рассматриваются отображения с прямыми и обратными модульными
условиями. В первой части статьи рассмотрены обратные модульные
неравенства, которым удовлетворяют прямые отображения, во второй
части -- напротив, модульные условия -- прямые, а рассматриваемые
отображения -- обратные. Круг изучаемых вопросов включает в себя
также глобальное поведение отображений (сходимость отображений в
замыкании заданной области). Основные определения и обозначения,
используемые в статье, могут быть найдены в монографиях~\cite{Ri} и
\cite{MRSY}, и потому опускаются.

Как известно, частью определения квазиконформных отображений
является неравенство
\begin{equation}\label{eq2}
M(\Gamma)\leqslant K\cdot M(f(\Gamma))\,,
\end{equation}
где $M$ -- модуль семейств кривых $\Gamma,$ а $K$ -- коэффициент
квазиконформности (см. \cite{Va}). Неравенство (\ref{eq2}) также
выполнено для квазиконформных отображений с ветвлением, которые
принято называть квазирегулярными отображениями (см.~\cite{Ri}). В
последнем случае, в (\ref{eq2}) следует взять $K:=N(f, A)K_O(f),$
где $K_O={\rm vrai\,} \sup K_O(x, f),$ $K_O(x, f)$ -- внешняя
дилатация, а $N(f, A)$ -- максимальная кратность отображения $f$ на
множестве $A,$ к которому принадлежат образы кривых семейства
$\Gamma$ (см. \cite[теорема~3.2]{MRV$_1$}, либо \cite[теорема~6.7,
гл.~II]{Ri}). По поводу отображений с ограниченным искажением
евклидового пространства также установлено, что они открыты и
дискретны и, в частности, они являются нульмерными отображениями,
см.~\cite[теорема~I.4.1]{Ri}. (Отображение $f:D\rightarrow
\overline{{\Bbb R}^n}={\Bbb R}^n\cup\{\infty\}$ называется {\it
нульмерным}, если ${\rm dim\,}\{f^{\,-1}(y)\}=0$ для каждого $y\in
\overline{{\Bbb R}^n},$ где ${\rm dim}$ обозначает топологическую
размерность множества, см. \cite{HW}).

\medskip
Рассмотрим следующий вопрос: будет ли произвольное отображение,
удовлетворяющее оценке (\ref{eq2}), открытым и дискретным
(нульмерным) ? Положительный ответ для пространства ${\Bbb R}^n$ дан
в сравнительно недавней работе \cite{Sev$_1$}, где рассматривается
даже более общее условие вида
\begin{equation}\label{eq2A}
M(\Gamma )\leqslant \int\limits_{D^{\,\prime}} Q(y)\cdot \rho_*^n
(y)\,d\mu^{\,\prime}(y)\,,
\end{equation}
где $Q$ -- заданная измеримая по Лебегу функция, удовлетворяющая
условию типа $FMO$ (см.~\cite{MRSY}), $\rho_*$ -- произвольная
неотрицательная борелевская функция такая, что
$\int\limits_{\gamma_*}\rho_*(y)ds\geqslant 1\quad\forall\quad
\gamma_*\in f(\Gamma)\,,$ область $D^{\,\prime}=f(D),$ а $M$ --
конформный модуль семейства кривых. Более того, достаточно, чтобы
$f$ было пределом последовательности отображений, удовлетворяющих
условию (\ref{eq2A}) (см. \cite[теорема]{Sev$_2$}), а вместо
показателя $n$ в правой части (\ref{eq2A}) можно брать произвольное
число $p\in (n-1, n].$

\medskip
В настоящей заметке мы хотим перенести указанный результат на
метрические пространства и, тем самым, подытожить наши исследования
в этом направлении. Ниже мы покажем, что аналогичное утверждение
имеет место в пространствах, регулярных по Альфорсу, в которых
выполнено так называемое $(1, p)$-неравенство Пуанкаре, $p\in
(\alpha-1, \alpha]$ и $\alpha$ -- хаусдорфова размерность
рассматриваемого метрического пространства. Вместо неравенства
(\ref{eq2A}) можно брать даже несколько более общее соотношение,
которое приведено в тексте ниже.

Всюду далее $\left(X,\,d, \mu\right)$ и
$\left(X^{\,\prime},\,d^{\,\prime}, \mu^{\,\prime}\right)$
--- произвольные  метрические пространства с метриками $d$ и $d^{\,\prime},$
наделённые локально конечными борелевскими мерами $\mu$ и
$\mu^{\,\prime},$ и конечными хаусдорфовыми размерностями
$\alpha\geqslant 2$ и $\alpha^{\,\prime}\geqslant 2,$
соответственно. Ниже мы считаем известными определения, связанные с
кривыми в метрическом пространстве, длинами дуг, интегралами,
условиями допустимости и так далее (см. \cite[разд.~13]{MRSY}).
Определения, связанные с  регулярностью по Альфорсу и $(1;
p)$-неравенство Пуанкаре могут быть найдены, напр., в~\cite[раздел
7.22]{He}.

\medskip
Пусть $G$ -- область в метрическом пространстве $(X,d,\mu)$. Следуя
\cite[разд.~13.4]{MRSY}, будем говорить, что локально интегрируемая
в $X$ функция $\varphi:G\rightarrow{\Bbb R}$ имеет {\it конечное
среднее колебание в точке $x_{0}\in\overline{G}$}, пишем $\varphi
\in FMO(x_{0})$, если $\overline{\lim\limits_{\varepsilon\rightarrow
0}}\,\, \,\frac{1}{\mu(G(x_{0},\varepsilon))}
\int\limits_{G(x_{0},\varepsilon)}|\varphi(x)-\overline{\varphi}_{\varepsilon}|\,\,d\mu(x)<\infty,$
где $\overline{\varphi}_{\varepsilon}
=\frac{1}{\mu(G(x_{0},\varepsilon))}
\int\limits_{G(x_{0},\varepsilon)}\varphi(x)\,\,d\mu(x)$ -- среднее
интегральное значение функции $\varphi(x)$ над множеством
$G(x_0,\varepsilon)=B(x_0, \varepsilon)\cap G=\{x\in X: d(x,
x_0)<\varepsilon\}\cap G$ по отношению к мере $\mu$. Пусть
$p\geqslant 1,$ тогда {\it $p$-модулем} семейства кривых $\Gamma$ в
метрическом пространстве $X$ называется величина
$$M_p(\Gamma)=\inf\limits_{\rho \in \,{\rm adm}\,\Gamma}
\int\limits_{X} \rho^p(x)\,d\mu(x)\,.$$
{\it Областью $D$ в метрическом пространстве $X$} называется
множество $D,$ являющееся линейно связным в $X.$ Пусть $E,$
$F\subset X$ -- произвольные множества. Обозначим через
$\Gamma(E,F,D)$ семейство всех кривых $\gamma:[a,b]\rightarrow X$
которые соединяют $E$ и $F$ в $D\,,$ т.е. $\gamma(a)\in
E,\gamma(b)\in\,F$ и $\gamma(t)\in D$ при $t \in (a, b).$ Пусть $D$
-- область в $X.$ Для $y_0\in f(D)$ и чисел $0<r_1<r_2<\infty$
обозначим
\begin{equation}\label{eq1**}
A(y_0, r_1, r_2)=\left\{ y\,\in\,X^{\,\prime}: r_1<d(y,
y_0)<r_2\right\}\,,\,\,\, S(y_0, r)=\{y\,\in\,X^{\,\prime}: d(y,
y_0)=r\}\,.\end{equation}
Пусть $G^{\,\prime}$ -- область в $X^{\,\prime}$ и
$Q:G^{\,\prime}\rightarrow [0, \infty]$ -- измеримая относительно
меры $\mu^{\,\prime}$ функция. Если $f:G\rightarrow G^{\,\prime}$ --
заданное отображение, то для фиксированного $y_0\in f(G)$ и
произвольных $0<r_1<r_2<\infty$ обозначим через $\Gamma(y_0, r_1,
r_2)$ семейство всех кривых $\gamma$ в области $G$ таких, что
$f(\gamma)\in \Gamma(S(y_0, r_1), S(y_0, r_2), A(y_0, r_1, r_2)).$
Для заданных $p, q> 1$ рассмотрим вместо (\ref{eq2}) неравенство
\begin{equation} \label{eq2*A}
M_p(\Gamma(y_0, r_1, r_2))\leqslant \int\limits_{G^{\,\prime}}
Q(y)\cdot \eta^q (d^{\,\prime}(y, y_0)) d\mu^{\,\prime}(y)\,,
\end{equation}
выполненное для любой неотрицательной измеримой по Лебегу функции
$\eta: (r_1,r_2)\rightarrow [0,\infty ]$ такой, что
\begin{equation}\label{eqA2}
\int\limits_{r_1}^{r_2}\eta(r) dr\geqslant 1\,.
\end{equation}
Пусть $(X,d)$ и $\left(X^{\,{\prime}},{d}^{\,{\prime}}\right)$ --
метрические пространства с расстояниями  $d$  и ${d}^{\,{\prime}},$
соответственно. Говорят, что последовательность отображений $f_k:X
\rightarrow X^{\,\prime}$, $k=1,2\dots$, сходится {\it локально
равномерно} к отображению $f:X\rightarrow {X}^{\,\prime}$, если
$\sup\limits_{x\,\in\,C}
d^{\,\prime}\left(f_k(x)\,,f(x)\right)\,\rightarrow\,0$
при $k\,\rightarrow\,\infty$ на любом компакте $C\subset X.$ Имеет
место следующая

\medskip
\begin{theorem}\label{th3}{\sl\, Предположим, $\left(X,\,d, \mu\right)$ и
$\left(X^{\,\prime},\,d^{\,\prime}, \mu^{\,\prime}\right)$ --
метрические пространства с метриками $d$ и $d^{\,\prime},$
наделённые локально конечными борелевскими мерами $\mu$ и
$\mu^{\,\prime},$ а $G$ и $G^{\,\prime}$ -- области в $X$ и
$X^{\,\prime},$ имеющие конечные хаусдорфовы размерности
$\alpha\geqslant 2$ и $\alpha^{\,\prime}\geqslant 2,$
соответственно. Предположим, кроме того, $G$ -- локально компактное
и локально связное пространство, $\alpha$-регулярное по Альфорсу, в
котором выполнено $(1; p)$-неравенство Пуанкаре при некотором $p\in
(\alpha-1, \alpha].$

Пусть $f_m:G\,\rightarrow\,{G}^{\,\prime},$ $m=1,2,\ldots$ --
последовательность непрерывных отображений, сходящаяся локально
равномерно к некоторому отображению
$f:G\,\rightarrow\,{G}^{\,\prime}.$ Тогда, если $f_m$ при каждом
$m\in {\Bbb N}$ удовлетворяет (\ref{eq2*A}) в каждой точке $y_0\in
f(G)$ при некотором $q\in (1, \alpha^{\,\prime}]$ и любой
неотрицательной измеримой по Лебегу функции $\eta:
(r_1,r_2)\rightarrow [0,\infty ],$ удовлетворяющей условию
(\ref{eqA2}), и $Q\in FMO$ в каждой точке $y_0\in f(G),$ то
отображение $f$ либо нульмерно, либо постоянно в $G.$
В частности, заключение теоремы~\ref{th3} имеет место, если
%
$$\limsup\limits_{\varepsilon\rightarrow
0}\frac{1}{\mu^{\,\prime}(B(y_0, \varepsilon))} \int\limits_{B(y_0,
\varepsilon)}Q(y)\,\,d\mu^{\,\prime}(y)<\infty\quad\forall\,\,y_0\in
f(G)\,.
$$
}
\end{theorem}

\medskip
{\bf 2. Формулировка и доказательство основной леммы. Доказательство
теоремы~\ref{th3}.} Связный компакт $C\subset X$ будем называть {\it
континуумом}. Следующая лемма включает в себя основной результат
настоящей работы в наиболее общей ситуации.

\medskip
\begin{lemma}\label{lem1}
{\sl\,Предположим, $\left(X,\,d, \mu\right)$ и
$\left(X^{\,\prime},\,d^{\,\prime}, \mu^{\,\prime}\right)$ --
метрические пространства с метриками $d$ и $d^{\,\prime},$
наделённые локально конечными борелевскими мерами $\mu$ и
$\mu^{\,\prime},$ а $G$ и $G^{\,\prime}$ -- области в $X$ и
$X^{\,\prime},$ имеющие конечные хаусдорфовы размерности
$\alpha\geqslant 2$ и $\alpha^{\,\prime}\geqslant 2,$
соответственно. Предположим, кроме того, $G$ -- локально компактное
и локально связное пространство,  $\alpha$-регулярное по Альфорсу, в
котором выполнено $(1; p)$-неравенство Пуанкаре при некотором $p\in
(\alpha-1, \alpha].$

Далее, предположим, что при некотором $q\in (1, \alpha^{\,\prime}]$
и каждого $y_0\in f(G)$ найдётся $\varepsilon_0>0,$ для которого
выполнено соотношение
\begin{equation} \label{eq4!}
\int\limits_{A(y_0, \varepsilon, \varepsilon_0)
}Q(y)\cdot\psi^q(d^{\,\prime}(y, y_0)) \
d\mu^{\,\prime}(y)\,=\,o\left(I^q(\varepsilon, \varepsilon_0)\right)
\end{equation}
для некоторой измеримой по Лебегу функции
$\psi(t):(0,\infty)\rightarrow (0,\infty),$ такой что
%
%
$$0< I(\varepsilon,
\varepsilon_0):=\int\limits_{\varepsilon}^{\varepsilon_0}\psi(t)\,dt
< \infty$$
%
%
при всех $\varepsilon\in(0,\varepsilon_0),$ где $A(y_0, \varepsilon,
\varepsilon_0)$ определено в (\ref{eq1**}) при $r_1=\varepsilon,$
$r_2=\varepsilon_0.$ Потребуем также, чтобы $I(\varepsilon,
\varepsilon_0)\rightarrow\infty$ при $\varepsilon\rightarrow 0.$

Пусть $f_m:G\,\rightarrow\,{G}^{\,\prime},$ $m=1,2,\ldots$ --
последовательность отображений, сходящаяся локально равномерно к
некоторому отображению $f:G\,\rightarrow\,{G}^{\,\prime}.$ Тогда,
если $f_m$ при каждом $m\in {\Bbb N}$ удовлетворяет (\ref{eq2*A}) в
каждой точке $y_0\in f(G)$ при некотором $q\in (1,
\alpha^{\,\prime}]$ и любой неотрицательной измеримой по Лебегу
функции $\eta: (r_1,r_2)\rightarrow [0,\infty ]$ с условием
(\ref{eqA2}), то $f$ либо нульмерно, либо постоянно.}
\end{lemma}

\medskip
\begin{proof}
Если отображение $f$ постоянно, доказывать нечего. Пусть
$f\not\equiv const.$ Предположим противное, а именно, что
отображение $f$ не нульмерно. Тогда найдётся $y_0\in G^{\,\prime},$
такое что множество $\{f^{\,-1}(y_0)\}$ не является всюду разрывным.
Следовательно, по определению, существует невырожденное связное
множество $C\subset \{f^{\,-1}(y_0)\}.$ Поскольку пространство $X$
локально компактно, можно считать, что $C$ -- континуум.

Поскольку по предположению $f\not\equiv y_0,$ ввиду непрерывности
отображения $f$ найдётся $x_0\in G$ и $\delta_0>0:$
$\overline{B(x_0, \delta_0)}\subset G$ и
%
$$f(x)\ne y_0\qquad 
\forall\quad x\in \overline{B(x_0, \delta_0)}\,.$$
%
Ввиду локальной компактности $G$ можно считать, что
$\overline{B(x_0, \delta_0)}$ -- компакт в $X.$ Кроме того, ввиду
локальной связности $G$ найдётся связная окрестность $U\subset
B(x_0, \delta_0).$ По определению, $U$ содержит некоторый шар
$B(x_0,\overline{\delta_0})\subset U.$ Заметим, что в силу
регулярности по Альфорсу метрического пространства $X$ шар
$B(x_0,\overline{\delta_0})$ не может быть одноточечным множеством.
Тогда $\overline{U}$ -- невырожденный континуум в $G.$

\medskip
В силу~\cite[предложение~4.7]{AS}, при $p\in (\alpha-1, \alpha]$
будем иметь:
\begin{equation}\label{eq4*}
M_p\left(\Gamma\left(C, \overline{U}, G\right)\right)>0\,.
\end{equation}
При достаточно большом $m\in {\Bbb N}$ рассмотрим семейство кривых
$f_m\left(\Gamma\left(C, \overline{U}, G\right)\right)\,.$
Заметим, что в силу локально равномерной сходимости $f_m$ к $f$
может быть построена подпоследовательность $f_{m_k}$ такая, что
$d^{\,\prime}(f_{m_k}(x), y_0)<1/2^k$ при всех $k\in {\Bbb N}$ и
всех $x\in C.$ С другой стороны, $f(\overline{U})$ -- компакт в
$X^{\,\prime}$ как непрерывный образ компакта, поэтому
$d^{\,\prime}\,(y_0, f(\overline{U}))\geqslant \sigma_0>0.$
Поскольку $f_m$ сходится к $f$ локально равномерно,
$$d^{\,\prime}(f_m(x), y_0)\geqslant d^{\,\prime}(f(x), y_0)-
d^{\,\prime}(f_m(x), f(x))\geqslant \sigma_0/2$$ при всех $x\in
\overline{U}$ и всех $m\geqslant m_0.$ В таком случае, каждая кривая
$\gamma\in f_{m_k}\left(\Gamma\left(C, \overline{U},
G\right)\right)$ имеет подкривую $\gamma^{\,\prime}\in \Gamma(S(y_0,
1/2^k), S(y_0, \sigma_0/2), A(y_0, 1/2^k, \sigma_0/2))$ при
достаточно больших $k\geqslant k_0$ (см.
\cite[предложение~13.3]{MRSY}). Отсюда $\Gamma\left(C, \overline{U},
G\right)>\Gamma_{f_{m_k}}(y_0, 1/2^k, \sigma_0/2)$ и, значит, ввиду
минорирования модуля (см.~\cite[теорема~1]{Fu})
\begin{equation}\label{eq1}
M_p\left(\Gamma\left(C, \overline{U}, G\right)\right)\leqslant
M_p(\Gamma_{f_{m_k}}(y_0, 1/2^k, \sigma_0/2))\,.
\end{equation}
Рассмотрим следующую функцию
$$\eta_{k}(t)=\left\{
\begin{array}{rr}
\psi(t)/I(1/2^k, \sigma_0/2), & t\in [1/2^k,
\sigma_0/2],\\
0,  & t\in {\Bbb R}\setminus [1/2^k, \sigma_0/2]\,,
\end{array}
\right. $$
где $I(1/2^k, \sigma_0/2)=\int\limits_{1/2^k}^{\sigma_0/2}\,\psi (t)
dt.$ Заметим, что функция $\eta_{k}$ удовлетворяет условию вида
(\ref{eqA2}) при $r_1=1/2^k$ и $r_2=\sigma_0/2.$ Тогда согласно
неравенствам (\ref{eq2*A}), (\ref{eq4!}) и (\ref{eq1}) мы получаем,
что
\begin{equation}\label{eq6}
M_p\left(\Gamma\left(C, \overline{U}, G\right)\right)\leqslant
\frac{1}{I^q(1/2^k, \sigma_0/2)}\int\limits_{A(y_0, 1/2^k,
\sigma_0/2)}Q(y)\psi^q(d^{\,\prime}(y,
y_0))\,d\mu^{\,\prime}(y)\,\leqslant
\end{equation}
$$\leqslant \frac{C}{I^q(1/2^k,
\varepsilon_0)}\int\limits_{A(y_0, 1/2^k,
\sigma_0/2)}Q(y)\psi^q(d^{\,\prime}(y, y_0))\,d\mu^{\,\prime}(y)
\rightarrow 0$$
при $k\rightarrow \infty.$
Однако, соотношение (\ref{eq6}) противоречит неравенству
(\ref{eq4*}). Полученное противоречие доказывает, что отображение
$f$ является нульмерным, что и требовалось доказать.
\end{proof}$\Box$

\medskip
{\it Доказательство теоремы~\ref{th3}} вытекает из леммы \ref{lem1}
и \cite[лемма~13.2]{MRSY}.

\medskip
{\bf 3. О сходимости обратных отображений в метрических
пространствах.} Теперь рассмотрим вопрос о сходимости отображений,
удовлетворяющих <<обратному>> к (\ref{eq2*A}) неравенству. Такие
результаты в пространстве ${\Bbb R}^n$ были получены в работе
\cite[теорема~6.1]{Sev$_3$}, однако, присутствующее здесь условие
фиксации двух точек области нас не вполне устраивает (среди
дробно-линейных автоморфизмов единичного круга на себя есть,
например, всего одно такое отображение). В настоящем тексте мы
отказываемся от упомянутого условия нормировки, заменяя его условием
${\rm diam\,}f(A)\geqslant \delta>0,$ которое будем требовать для
всех отображений $f$ из рассматриваемого класса и фиксированного
континуума $A.$ Предположим, $f:D\rightarrow D^{\,\prime}$ --
фиксированное отображение, $\alpha$ и $\alpha^{\,\prime}$ --
хаусдорфовы размерности областей $D\subset X$ и $D^{\,\prime}\subset
X^{\,\prime},$ соответственно, и пусть теперь вместо требования
(\ref{eq2*A}) выполнено более сильное условие
\begin{equation}\label{eq1C}
M_{\alpha^{\,\prime}}(f(\Gamma))\leqslant \int\limits_DQ(x)\cdot
\rho^{\alpha}(x)\,d\mu(x)\,,
\end{equation}
где $Q:D\rightarrow [1, \infty]$ -- фиксированная измеримая функция,
а $\rho: D\rightarrow [0,\infty ]$ пробегает класс борелевских
функций, подчинённых неравенству
$ \int\limits_{\gamma}\rho(x)|dx|\geqslant 1.$
Будем говорить, что $f$ -- {\it $Q$-гомеоморфизм} в $D,$ если $f$
удовлетворяет условию (\ref{eq1C}) для каждого семейства кривых
$\Gamma$ в $D$ и произвольной $\rho\in{\rm adm}\,\Gamma.$

\medskip
В сравнительно недавних работах \cite{RSa} и \cite{Sm} был решён
вопрос о возможности непрерывного продолжения отображений данного
вида в метрических пространствах (см. \cite[лемма~6.1 и
теорема~6.1]{RSa} и \cite[лемма~5 и теорема~3]{Sm}). Также в работе
первого автора получено свойство равностепенной непрерывности
указанных отображений в замыкании области в ${\Bbb R}^n$ (см.
\cite[теорема~6.1]{Sev$_3$}). Целью настоящего раздела является
описание сходимости отображений с условием (\ref{eq1C}) в
метрических пространствах.

Напомним некоторые определения. Пусть $(X, d, \mu)$ -- метрическое
пространство с мерой $\mu.$ Определим {\it функцию Лёвнера
$\phi_{\alpha}:(0, \infty)\rightarrow [0, \infty)$ на $X$} по
следующему правилу:
\begin{equation}\label{eq2H}
\phi_{\alpha}(t)=\inf\{M_{\alpha}(\Gamma(E, F, X)): \Delta(E,
F)\leqslant t\}\,,
\end{equation}
где $\inf$ берётся по всем произвольным невырожденным
непересекающимся континуумам $E, F$ в $X,$ относительно которых
величина $\Delta(E, F)$ определяется так:
\begin{equation}\label{eq1H}
\Delta(E, F):=\frac{{\rm dist}\,(E, F)}{\min\{{\rm diam\,}E, {\rm
diam\,}F\}}\,.
\end{equation}
Пространство $X$ называется {\it пространством Лёвнера,} если
функция $\phi_n(t)$ положительна при всех положительных значениях
$t$ (см. \cite[разд.~2.5]{MRSY} либо \cite[гл.~8]{He}). Область $D$
в $X$ будем называть {\it областью квазиэкстремальной длины
относительно $p$-модуля}, сокр. {\it $QED$-облас\-тью}, если
$M_{\alpha}(\Gamma(E, F, X))\leqslant A\cdot M_{\alpha}(\Gamma(E, F,
D))$
для конечного числа $A\geqslant 1$ и всех континуумов $E$ и $F$ в
$D.$ Область $D$ будет называться {\it локально линейно связной в
точ\-ке} $x_0\in \overline{D},$ если для любой окрестности $U$ точки
$x_0$ найдется окрестность $V\subset U$ такая, что множество $V\cap
D$ линейно связно. В частности, будем говорить, что $D$ локально
линейно связна на границе $\partial D,$ если $D$ локально линейно
связна в каждой точке $x_0\in \partial D.$ Докажем следующее весьма
важное утверждение.

\medskip
\begin{lemma}\label{lem2}
{\sl\, Предположим, что $D\subset X$ и $D^{\,\prime}\subset
X^{\,\prime}$ -- области с конечными хаусдорфовыми размерностями
$\alpha\geqslant 2$ и $\alpha^{\,\prime}\geqslant 2,$
соответственно, $X^{\,\prime}$ -- $\alpha^{\,\prime}$-регулярное по
Альфорсу пространство Лёвнера. Пусть также $D^{\,\prime}$ является
$QED$-об\-лас\-тью. Тогда:

\medskip
{\textbf 1)} в $D^{\,\prime}$ имеет место свойство сближающихся
континуумов: если $E_k, F_k$ -- произвольные континуумы в
$D^{\,\prime}$ такие, что $\min\{{\rm diam}\, E_k, {\rm diam}\,
F_k\}\geqslant \delta,$ где $\delta>0$ -- фиксированное число, и
${\rm dist\,}(E_k, F_k)\rightarrow 0$ при $k\rightarrow\infty,$ то
$M_{\alpha^{\,\prime}}(\Gamma(E_k, F_k, D^{\,\prime}))\rightarrow
\infty$ при $k\rightarrow\infty;$

{\textbf 2)} граница области $D^{\,\prime}$ является слабо плоской,
т.е. какова бы ни была точка $x_0\in
\partial D^{\,\prime},$ для каждого $P>0$ и для любой окрестности $U$ точки
$x_0$ найдётся окрестность $V\subset U$ этой же точки такая, что
$M_{\alpha^{\,\prime}}(\Gamma(E, F, D^{\,\prime}))>P$ для
произвольных континуумов $E, F\subset D^{\,\prime},$ пересекающих
$\partial U$ и $\partial V.$

\medskip
Предположим, кроме того, что область $D$ локально линейно связна на
$\overline{D},$ тогда

\medskip
{\textbf 3)} если $U$ -- окрестность континуума $E_0\subset
{\overline{D}},$ то найдётся окрестность $V\subset U$ континуума
$E_0$ такая, что $V\cap D$ -- линейно связное множество.

\medskip
Пусть, кроме того, $\overline{D}$ и $\overline{D^{\,\prime}}$ --
компакты в $X$ и $X^{\,\prime},$ соответственно, и $Q\in L^1(D).$
Тогда:

\medskip
{\textbf 4)}  если $f:D\rightarrow D^{\,\prime}$ -- $Q$-гомеоморфизм
области $D$ на область $D^{\,\prime},$ то $g=f^{\,-1}$ продолжается
до непрерывного отображения $\overline{g}:\overline{D}\rightarrow
\overline{D^{\,\prime}},$ при этом,
$\overline{g}(\overline{D^{\,\prime}})=\overline{D}.$

{\textbf 5)} если никакая связная компонента границы $\partial
D^{\,\prime}$ не вырождается в точку и $f_m:D\rightarrow
D^{\,\prime}$ -- последовательность $Q$-гомеоморфизмов области $D$
на область $D^{\,\prime},$ удовлетворяющих для некоторого
(фиксированного) континуума $A\subset D$ условию ${\rm diam\,}
f_m(A)\geqslant \delta>0$ при всех $m=1,2,\ldots ,$ то найдётся
$\delta_1>0$ такое, что ${\rm dist\,}(f_m(A),
\partial D^{\,\prime})>\delta_1>0$ для всех $m\in {\Bbb N}.$ }
\end{lemma}

\medskip
\begin{proof}
Установим вначале свойство {\textbf 1).} Поскольку $X^{\,\prime}$ по
предположению является пространством Лёвнера и, кроме того,
$X^{\,\prime}$ является $\alpha^{\,\prime}$-регулярным по Альфорсу,
мы получим: $\phi_{\alpha^{\,\prime}}(t)\rightarrow \infty$ при
$t\rightarrow 0$ (см. \cite[теорема 8.23]{He}). Возьмём произвольное
$\varepsilon>0$ и для него найдём $t_0=t_0(\varepsilon)$ такое, что
при $t\in (0, t_0)$ выполнено
$\phi_{\alpha^{\,\prime}}(t)>\varepsilon.$ Положим $\Delta(E_k,
F_k)=t,$ где $\Delta(E_k, F_k)$ определено в (\ref{eq1H}). Тогда
ввиду (\ref{eq2H})
\begin{equation}\label{eq3H}
\varepsilon<M_{\alpha^{\,\prime}}(\Gamma(E_k, F_k, X^{\,\prime}))
\end{equation}
как только $t\in (0, t_0)$ и $\Delta(E_k, F_k)=t.$ Заметим, что
$\Delta(E_k, F_k)\leqslant\frac{1}{\delta}{\rm dist}\, (E_k, F_k),$
и если ${\rm dist}\, (E_k, F_k)\in (0, t_0\delta),$ то $\Delta(E_k,
F_k)=t\in (0, t_0)$ и, значит, имеет место соотношение (\ref{eq3H}).
Окончательно, для произвольного $\varepsilon>0$ нашлось
$t^{\,\prime}_0=t_0\delta$ такое, что как только ${\rm dist\,}(E_k,
F_k)\in (0, t^{\,\prime}),$ выполняется условие (\ref{eq3H}). Из
того, что $D^{\,\prime}$ является $QED$-областью (либо,
соответственно, $QED$-областью относительно
$\overline{D^{\,\prime}}$), вытекает, что
%
$$\varepsilon /L<M_{\alpha^{\,\prime}}(\Gamma(E_k, F_k,
D^{\,\prime}))\,,$$
%
где $L$ -- некоторая фиксированная постоянная, откуда и следует, что
$M_{\alpha^{\,\prime}}(\Gamma(E_k, F_k, D^{\,\prime}))\rightarrow
\infty$ при ${\rm dist\,}(E_k, F_k)\rightarrow 0,$
$k\rightarrow\infty.$

\medskip
Теперь докажем свойство {\textbf 2)}. Воспользуемся установленным
уже свойством ${\textbf 1)}.$ Пусть теперь $x_0\in\partial
D^{\,\prime}$ -- произвольная точка. Берём произвольную окрестность
$U$ точки $x_0$ и произвольное $P>0.$ Для числа $k\in {\Bbb N}$
найдём окрестность $V_k$ точки $x_0,$ лежащую в шаре
$\overline{B(x_0, 2^{-k})}.$ Рассмотрим континуумы $E$ и $F,$
пересекающие $\partial U$ и $\partial V_k.$ Заметим, что $\min\{{\rm
diam\, E}, {\rm diam}\, F\}\geqslant \delta>0$ при достаточно
больших $k,$ поскольку $U$ -- фиксированная окрестность, ${\rm
dist}\,(\partial V_k, x_0)\rightarrow 0 $ при $k\rightarrow \infty,$
а ${\rm diam\, E}$ и ${\rm diam}\, F$ не меньше расстояния между
$\partial U$ и $\partial V_k.$ Кроме того, ${\rm dist}\, (E,
F)\leqslant {\rm diam}\,\partial V_k\rightarrow 0$ при $k\rightarrow
\infty.$ Тогда в силу свойства ${\textbf 1)}$ имеем:
$M_{\alpha^{\,\prime}}(\Gamma(E, F,
D^{\,\prime}))=\alpha_k\rightarrow\infty$ при $k\rightarrow \infty.$
Подберём $k_0$ так, чтобы $\alpha_k>P$ при $k\geqslant k_0$ (это
число $k_0$ полностью определяется числом $P$). Положим
$V:=V_{k_0}.$ Тогда получаем, что $M_{\alpha^{\,\prime}}(\Gamma(E,
F, D^{\,\prime}))>P$ для произвольных континуумов $E, F\subset
D^{\,\prime},$ пересекающих $\partial U$ и $\partial V,$ что и
требовалось установить.

\medskip
Чтобы установить свойство {\textbf 3)}, будем следовать схеме
доказательства~\cite[лемма~2.2]{HK$_1$}. Поскольку $E_0$ --
континуум, а $D$ -- локально линейно связна на $\overline{D},$
найдётся конечно число номеров $1, 2, \ldots m$ и соответствующие
окрестности $V_1, V_2,\ldots V_m\subset U$ некоторых точек
$x_1,\ldots x_m\in E_0$ такие, что $W_j:=V_j\cap D$ -- связное
множество при каждом $j=1,2,\ldots, m,$ $V:=V_1\cup
V_2\cup\ldots\cup V_m\subset U$ и $V$ -- окрестность $E_0.$ Не
ограничивая общности, можно считать, что $V_i$ -- открытые множества
при каждом $i=1,2,\ldots m.$ Покажем, что $W:=V\cap D$ -- связное
множество. Предположим противное, тогда ввиду \cite[теорема~2.I,
гл.~5, \S\, 46]{Ku} найдутся два открытых непересекающихся множества
$G$ и $H\subset X$ такие, что $W=G\cup H,$ $G\cap W\ne\varnothing\ne
H\cap W.$ Не ограничивая общности, можно считать, что
$G=\bigcup\limits_{i=1}^kW_i$ и $H=\bigcup\limits_{i=k+1}^mW_i,$
$1\leqslant k<m.$ Поскольку $E_0\subset \overline{W}\subset
\overline{G}\cup\overline{H}$ и $E_0$ -- континуум,
$E_0\cap\overline{G}\cup\overline{H}\ne\varnothing,$ что вытекает из
определения связности для множества $E_0$ (см. \cite[I.5.46]{Ku}). В
таком случае, найдутся $1\leqslant i\leqslant k$ и $k<j\leqslant m$
такие, что $E_0\cap\overline{W_i}\cap\overline{W_j}\ne\varnothing.$
Пусть $x_0\in E_0\cap\overline{W_i}\cap\overline{W_j}.$ Поскольку
$V$ -- окрестность $E_0,$ найдётся $p\in {\Bbb N},$ $1\leqslant
p\leqslant m,$ и число $\varepsilon_0>0$ такие, что $B(x_0,
\varepsilon_0)\subset V_p.$ Поскольку $x_0\in
E_0\cap\overline{W_i},$ мы имеем $z_1\in W_i\cap B(x_0,
\varepsilon_0)\subset W_i\cap V_p\subset W_i\cap W_p.$ Аналогично,
найдётся $z\in W_j\cap W_p,$ следовательно, $W_i\cap
W_p\ne\varnothing\ne W_j\cap W_p.$ Значит, $G\cap H\ne\varnothing,$
что противоречит полученному выше. Полученное противоречие указывает
на связность множества $W:=V\cap D.$ В таком случае, $W$ -- линейно
связно ввиду~\cite[предложение~13.1]{MRSY}.

\medskip Свойство {\textbf 4)}, за исключением равенства
$\overline{g}(\overline{D^{\,\prime}})=\overline{D},$
$g:=f^{\,-1}=\overline{g}|_{D^{\,\prime}},$ вытекает из пункта
{\textbf 2)} и \cite[теорема~3]{Sm}. Легко видеть, что указанное
равенство $\overline{g}(\overline{D^{\,\prime}})=\overline{D}$ также
выполняется. В самом деле, по определению, мы имеем включение
$\overline{g}(\overline{D^{\,\prime}})\subset\overline{D},$ так что
для завершения доказательства нам достаточно установить обратное
включение $\overline{D}\subset
\overline{g}(\overline{D^{\,\prime}}).$ Если $x_0\in \overline{D},$
то для него найдётся $y_0\in D:$ $x_0=\overline{g}(y_0),$ поскольку
по условию $f(D)=D^{\,\prime}.$ Пусть теперь $x_0\in
\partial D,$ тогда найдутся последовательности $x_k\in D$
и $y_k\in D^{\,\prime}$ такие, что $x_k=\overline{g}(y_k)$ и
$x_k\rightarrow x_0$ при $k\rightarrow\infty.$ Поскольку
$\overline{D^{\,\prime}}$ -- компакт, то можно считать, что
$y_k\rightarrow y_0\in \overline{D^{\,\prime}}$ при
$k\rightarrow\infty.$ Так как $f:=g^{\,-1}$ -- гомеоморфизм, то
$y_0\in
\partial D^{\,\prime}.$ Поскольку $\overline{g}^{\,-1}$ непрерывно в $\overline{D^{\,\prime}},$ мы получим,
что $\overline{g}(y_k)\rightarrow g(y_0).$ Однако, в таком случае,
$\overline{g}(y_0)=x_0,$ ибо $\overline{g}(y_k)=x_k$ и
$x_k\rightarrow x_0,$ $k\rightarrow\infty.$ Значит,
$\overline{D}\subset f(\overline{D^{\,\prime}}),$ ввиду чего
$\overline{D}=\overline{g}(\overline{D^{\,\prime}}),$ что и
требовалось установить.

\medskip
Установим, наконец, свойство {\textbf 5)}. Предположим противное,
т.е., предположим что для каждого $k\in {\Bbb N}$ существует
$m=m_k:$ ${\rm dist\,}(f_{m_k}(A),
\partial D^{\,\prime})<1/k.$ Без ограничения общности
мы можем считать последовательность $m_k$ монотонно возрастающей. По
условию $\overline{D^{\,\prime}}$ -- компакт, поэтому и $\partial
D^{\,\prime}$ также компакт как замкнутое подмножество компакта
$\overline{D^{\,\prime}}.$ Кроме того, $f_{m_k}(A)$ компакт как
непрерывный образ компакта $A$ при отображении $f_{m_k}.$ Тогда
найдутся $x_k\in f_{m_k}(A)$ и $y_k\in\partial D^{\,\prime}$ такие,
что ${\rm dist\,}(f_{m_k}(A),
\partial D^{\,\prime})=d^{\,\prime}(x_k, y_k)<1/k$ (см. рисунок~\ref{fig1}).
\begin{figure}[h]
\centerline{\includegraphics[scale=0.6]{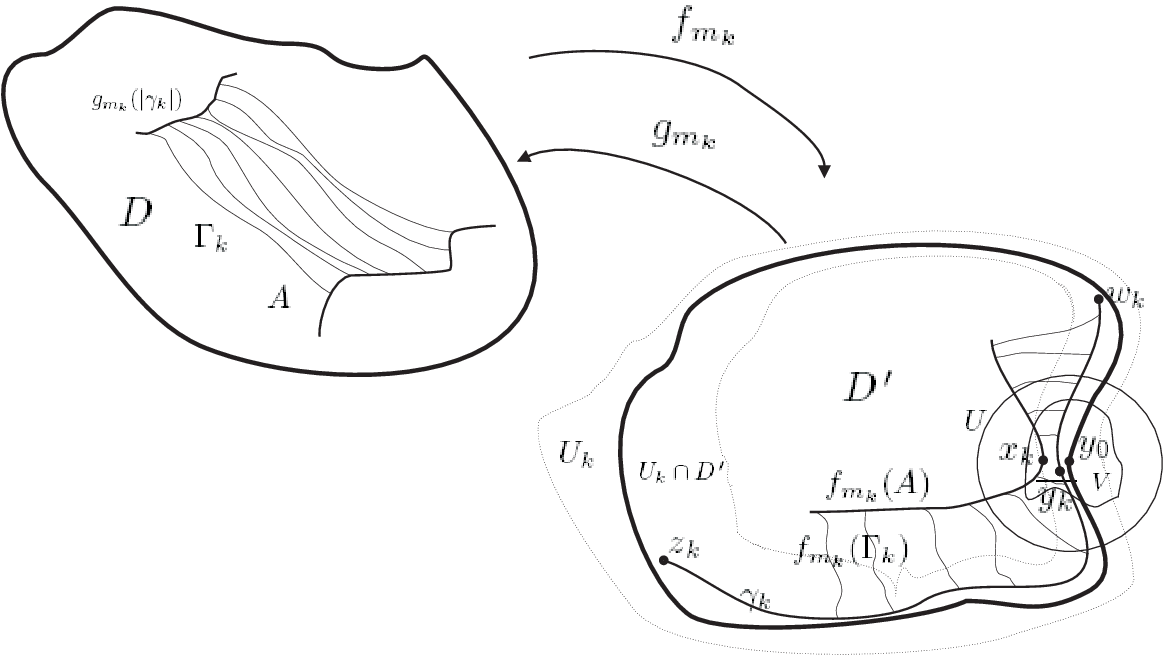}} \caption{К
доказательству пункта~\textbf{5)} леммы~\ref{lem2}.}\label{fig1}
\end{figure}
Так как $\partial D^{\,\prime}$ -- компакт, можно считать, что
$y_k\rightarrow y_0\in \partial D^{\,\prime},$ $k\rightarrow
\infty.$ Пусть $K_0$ -- связная компонента $\partial D^{\,\prime},$
содержащая точку $y_0.$ По условию $K_0$ -- невырожденный континуум
в $\partial D^{\,\prime},$ так что ${\rm diam}\,K_0>a_0>0.$

Согласно пункту~\textbf{4)}, при каждом $k\in {\Bbb N}$ отображение
$g_{m_k}:=f_{m_k}^{\,-1}$ продолжатся до непрерывного отображения
$\overline{g}_{m_k}:\overline{D^{\,\prime}}\rightarrow
\overline{D},$ более того, $\overline{g}_{m_k}$ равномерно
непрерывно на $\overline{D^{\,\prime}}.$ Тогда для всякого
$\varepsilon>0$ найдётся $\delta_k=\delta_k(\varepsilon)<1/k$ такое,
что
\begin{equation}\label{eq3B}
d(\overline{g}_{m_k}(x), \overline{g}_{m_k}(x_0))<\varepsilon \quad
\forall\,\, x,x_0\in \overline{D^{\,\prime}},\quad d^{\,\prime}(x,
x_0)<\delta_k\,, \quad \delta_k<1/k\,.
\end{equation}
Пусть далее $\varepsilon>0$ -- произвольное число с условием
\begin{equation}\label{eq5B}
\varepsilon<(1/2)\cdot {\rm  dist}\,(\partial D, A)\,,
\end{equation}
где $A$ -- континуум из условия леммы. При каждом фиксированном
$k\in {\Bbb N}$ рассмотрим множество
$B_k:=\bigcup\limits_{x_0\in K_0}B(x_0, \delta_k)\,,\,k\in {\Bbb
N}.$
Заметим, что $B_k$ -- открытое множество, содержащее $K_0,$ другими
словами, $B_k$ -- некоторая окрестность континуума $K_0.$ Ввиду
пункта~\textbf{3)}, существует окрестность $U_k\subset B_k$
континуума $K_0,$ такая, что $U_k\cap D^{\,\prime}$ линейно связно.
Пусть ${\rm diam}\,K_0=m_0,$ тогда найдутся $z_0, w_0\in K_0$ такие,
что ${\rm diam}\,K_0=d^{\,\prime}(z_0, w_0)=m_0.$ Следовательно,
можно выбрать последовательности $\overline{y_k}\in U_k\cap
D^{\,\prime},$ $z_k\in U_k\cap D^{\,\prime}$ и $w_k\in U_k\cap
D^{\,\prime}$ так, что $z_k\rightarrow z_0,$
$\overline{y_k}\rightarrow y_0$ и $w_k\rightarrow w_0$ при
$k\rightarrow\infty.$ Можно считать, что
\begin{equation}\label{eq2C}
d(z_k, w_k)>m_0/2,\quad \forall\,\, k\in {\Bbb N}\,.
\end{equation}
Соединим последовательно точки $z_k,$ $\overline{y_k}$ и $w_k$
кривой $\gamma_k$ в $U_k\cap D^{\,\prime}$ (это возможно, поскольку
$U_k\cap D^{\,\prime}$ линейно связно). Пусть $|\gamma_k|$ -- как
обычно, носитель (образ) кривой $\gamma_k$ в $D^{\,\prime}.$ Тогда
$g_{m_k}(|\gamma_k|)$ -- компакт в $D.$ Пусть $x\in|\gamma_k|,$
тогда найдётся $x_0\in K_0:$ $x\in B(x_0, \delta_k).$ Зафиксируем
$\omega\in A\subset D.$ Поскольку $x\in|\gamma_k|,$ то $x$ --
внутренняя точка области $D^{\,\prime},$ так что мы вправе писать
$g_{m_k}(x)$ вместо $\overline{g}_{m_k}(x)$ для указанных $x.$ Из
(\ref{eq3B}) и (\ref{eq5B}), ввиду неравенства треугольника для
достаточно больших $k\in {\Bbb N},$ получаем:
$$d(g_{m_k}(x),\omega)\geqslant
d(\omega,
\overline{g}_{m_k}(x_0))-d(\overline{g}_{m_k}(x_0),g_{m_k}(x))\geqslant$$
\begin{equation}\label{eq4}
\geqslant {\rm  dist}\,(\partial D, A)-(1/2)\cdot{\rm
dist}\,(\partial D, A)={\rm  dist}\,(\partial D, A)>\varepsilon\,.
\end{equation}
Переходя в (\ref{eq4}) к $\inf$ по всем $x\in |\gamma_k|$ и всем
$\omega\in A,$ мы получим:
\begin{equation}\label{eq6B}
{\rm dist}\,(g_{m_k}(|\gamma_k|), A)>\varepsilon, \quad\forall\,\,
k=1,2,\ldots \,.
\end{equation}
Ввиду (\ref{eq6B}) длина произвольной кривой, соединяющей компакты
$g_{m_k}(|\gamma_k|)$ и $A$ в $D,$ не меньше $\varepsilon.$ Положим
$\Gamma_k:=\Gamma(g_{m_k}(|\gamma_k|), A, D),$ тогда функция
$\rho(x)=1/\varepsilon$ при $x\in D$ и $\rho(x)=0$ при $x\not\in D$
допустима для $\Gamma_k.$ По определению отображений $f_{m_k}$ в
(\ref{eq1C}) имеем:
\begin{equation}\label{eq4B}
M_{\alpha^{\,\prime}}(f_{m_k}(\Gamma_k))\leqslant
\frac{1}{\varepsilon^{\,\alpha}}\int\limits_DQ(x)\,d\mu(x)=c=c(\varepsilon,
Q)<\infty\,,
\end{equation}
поскольку по условию $Q\in L^1(D).$
Однако, соотношение~(\ref{eq4B}) противоречит пункту~\textbf{1)}. В
самом деле, $\Gamma(f_{m_k}(A), |\gamma_k|,
D^{\,\prime})=f_{m_k}(\Gamma(A, g_{m_k}(|\gamma_k|),
D))=f_{m_k}(\Gamma_k),$ ${\rm diam}\,f_{m_k}(A)\geqslant \delta$ по
условию, ${\rm diam}\,|\gamma_k|\geqslant d(z_k, w_k)>m_0/2$ ввиду
(\ref{eq2C}), кроме того, ${\rm dist}\,(f_{m_k}(A),
|\gamma_k|)\leqslant d^{\,\prime}(x_k, \overline{y_k})\rightarrow 0$
при $k\rightarrow\infty,$ поскольку каждая из последовательностей
$x_k$ и $\overline{y_k}$ сходятся при $k\rightarrow\infty$ к точке
$y_0.$ В силу пункта~\textbf{1)} тогда
$M_{\alpha^{\,\prime}}(\Gamma(f_{m_k}(A), |\gamma_k|,
D^{\,\prime}))=M_{\alpha^{\,\prime}}(f_{m_k}(\Gamma_k))\rightarrow\infty$
при $k\rightarrow\infty,$ что противоречит~(\ref{eq4B}). Полученное
противоречие опровергает предположение о неравенстве ${\rm
dist\,}(f_{m_k}(A),
\partial D^{\,\prime})<1/k.$ Лемма доказана.
~$\Box$
\end{proof}

\medskip
Будем рассматривать в дальнейшем области $D\subset X$,
удовлетворяющие следующему условию~$\textbf{A}:$ {\it любые две пары
точек $a\in D, b\in \overline{D},$ и $c\in D, d\in \overline{D}$
можно соединить непересекающимися между собой кривыми $C_1$ и $C_2$
в области $D.$} В настоящей работе мы покажем, что области в ${\Bbb
R}^n,$ $n\geqslant 2,$ с локально связной границей всегда
удовлетворяют условию~$\textbf{A}$ (см. предложение~\ref{pr1}).

\medskip Семейство $\frak{F}$ отображений $f:X\rightarrow {X}^{\,\prime}$
называется {\it равностепенно непрерывным в точке} $x_0 \in X,$ если
для любого $\varepsilon>0$ найдётся $\delta>0$ такое, что
${d}^{\,\prime} \left(f(x),f(x_0)\right)<\varepsilon$ для всех $x$
таких, что $d(x,x_0)<\delta$ и для всех $f\in \frak{F}.$ Семейство
$\frak{F}$ {\it равностепенно непрерывно}, если $\frak{F}$
равностепенно непрерывно в каждой  точке $x_0\in X.$ Для областей
$D\subset X,$ $D^{\,\prime}\subset X^{\,\prime}$ и произвольной
измеримой по Лебегу функции $Q: X\rightarrow [1, \infty],$
$Q(x)\equiv 0$ при $x\not\in D,$ обозначим через ${\frak R}_Q(D,
D^{\,\prime})$ семейство всех гомеоморфизмов
$g:D^{\,\prime}\rightarrow D$ области $D^{\,\prime}$ на область $D$
таких, что $f=g^{\,-1},$ $f:D\rightarrow D^{\,\prime}$ --
$Q$-гомеоморфизм в $D.$ Справедливо следующее утверждение.

\medskip
\begin{theorem}\label{th4A} {\sl Предположим, что $D\subset X$ и $D^{\,\prime}\subset X^{\,\prime}$
-- области с конечными хаусдорфовыми размерностями $\alpha\geqslant
2$ и $\alpha^{\,\prime}\geqslant 2,$ соответственно. Пусть также:

1) пространство $X^{\,\prime}$ является
$\alpha^{\,\prime}$-регулярным по Альфорсу пространством Лёвнера,

2) область $D$ локально линейно связна на $\overline{D},$
$\overline{D}$ и $\overline{D^{\,\prime}}$ -- компакты в $X$ и
$X^{\,\prime},$ соответственно, кроме того, $\partial D$ содержит не
менее двух точек,

3) область $D^{\,\prime}$ является $QED$-об\-лас\-тью,

4) выполнено условие~\textbf{A},

5) $Q\in L^1(D).$

\medskip
Тогда семейство ${\frak R}_Q(D, D^{\,\prime})$ является
равностепенно непрерывным в $D^{\,\prime}.$}
\end{theorem}

\medskip
\begin{proof}
Проведём доказательство теоремы~\ref{th4A} от противного.
Предположим, что семейство ${\frak R}_Q(D, D^{\,\prime})$ не
является равностепенно непрерывным в некоторой точке $y_0\in
D^{\,\prime},$ другими словами, найдутся $y_0\in D^{\,\prime}$ и
$\varepsilon_0>0,$ такие что для любого $m\in {\Bbb N}$ существует
элемент $y_m\in D^{\,\prime}$ с условием $d^{\,\prime}(y_m,
y_0)<1/m$ и гомеоморфизм $g_m\in{\frak R}_Q(D, D^{\,\prime})$ такие
что
\begin{equation}\label{eq13***}
d(g_m(y_m), g_m(y_0))\geqslant \varepsilon_0\,.
\end{equation}
Поскольку по условию $\overline{D}$ является компактом, мы можем
считать, что последовательности $g_m(y_m)$ и $g_m(y_0)$ сходятся при
$m\rightarrow\infty$ к точкам $\overline{x_1}$ и $\overline{x_2}\in
\overline{D}.$ В силу неравенства (\ref{eq13***}) по непрерывности
метрики $d(\overline{x_1}, \overline{x_2})\geqslant \varepsilon_0.$
Соединим точку $\overline{x_1}$ с точкой $x_1\in \partial D,$ а
точку $\overline{x_2}$ -- с точкой $x_2\in \partial D,$ $x_1\ne
x_2,$ непересекающимися кривыми $\gamma_1$ и $\gamma_2,$
соответственно, так что $\gamma_i(t)\in D$ при $0<t<1,$
$\gamma_1(0)=\overline{x_1},$ $\gamma_1(1)=x_1,$
$\gamma_2(0)=\overline{x_2},$ $\gamma_2(1)=x_2$ (это возможно по
условию теоремы; см. рисунок~\ref{fig2}). В случае, если одна из
точек $\overline{x_1}$ или $\overline{x_2}$ граничная, отрезки
$\gamma_1$ и $\gamma_2$ вырождаются в точку по определению.

Пусть $U_1$ и $U_2$ -- непересекающиеся окрестности точек
$\overline{x_i},$ $i=1,2,$ такие, что $W_i:=U_i\cap D$ -- связное
множество; пусть также $P_i$ -- непересекающиеся окрестности точек
$x_i,$ $i=1,2,$ такие что $L_i:=P_i\cap D$ -- связное множество (все
такие окрестности существуют, поскольку по условию $D$ локально
связна на $\partial D$). Если $\overline{x_i},$ $i=1,2,$ --
внутренние точки области $D,$ то можно выбрать окрестности так, что
$\overline{L_i}\cap \overline{W_j}=\varnothing,$ $i,j=1,2.$ Положим
$f_m:=g_m^{\,-1}.$ Построим две последовательности $z^1_m\rightarrow
x_1$ и $z^2_m\rightarrow x_2$ при $m\rightarrow\infty,$
$m=1,2,\ldots ,$ следующий образом. Так как $f_m$ -- гомеоморфизм,
то предельное множество $C(f_m, x_1)$ лежит на $\partial
D^{\,\prime}$ (см. \cite[предложение~13.5]{MRSY}), поэтому найдётся
точка $z^1_m\in D$ такая, что ${\rm dist}\,(f_m(z^1_m),
\partial D^{\,\prime})<1/m.$ Так как $\overline{D^{\,\prime}}$ --
компакт, то можно считать, что последовательность
$f_m(z^1_m)\rightarrow p_1\in \partial D^{\,\prime}$ при
$m\rightarrow\infty.$ Аналогично строим последовательность $z^2_m:$
можно считать, что $f_m(z^2_m)\rightarrow p_2\in \partial
D^{\,\prime}$ при $m\rightarrow\infty.$ Можно считать, что $z^1_m\in
L_1$ и $z^2_m\in L_2$ при всех $m\in {\Bbb N.}$

Если точка $\overline{x_1}\in D,$ то последовательность $z^1_m\in D$
можно выбрать такую, что $z^1_m\in D\in |\gamma_1|;$ в этом случае,
пусть $P^1_m$ обозначает часть кривой $|\gamma_1|,$ соединяющая
точки $\overline{x_1}$ и $z^1_m.$ Положим тогда $P_m:=P^1_m\cup
\overline{W_1}.$ В противном случае, если $\overline{x_1}\in
\partial D,$ то пусть $P_m$ -- кривая, соединяющая точки $z^1_m$ и
$g_m(y_m)$ в $W_1.$

Аналогично, если точка $\overline{x_2}\in D,$ то последовательность
$z^2_m\in D$ можно выбрать такую, что $z^2_m\in D\in |\gamma_2|;$ в
этом случае, пусть $P^2_m$ обозначает часть отрезка $|\gamma_2|,$
соединяющая точки $\overline{x_2}$ и $z^2_m.$ Положим тогда
$Q_m:=P^2_m\cup \overline{W_2}.$ В противном случае, если
$\overline{x_2}\in
\partial D,$ то пусть $Q_m$ -- кривая, соединяющая точки $z^2_m$ и
$g_m(y_0)$ в $W_2.$
\begin{figure}[h]
\centerline{\includegraphics[scale=0.6]{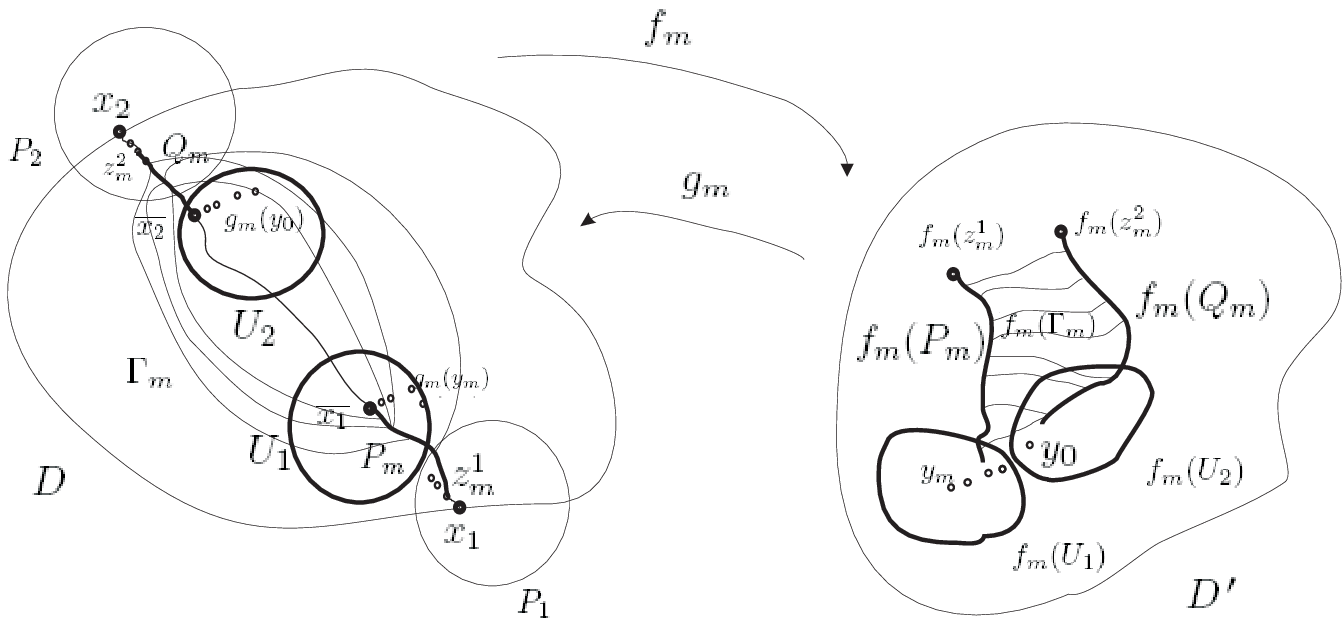}} \caption{К
доказательству теоремы~\ref{th4A}}\label{fig2}
\end{figure}
По построению, $P_m$ и $Q_m$ -- непересекающиеся континуумы в $D$
такие, что $l_m:={\rm dist}\,(P_m, Q_m)>l>0.$ Пусть
$\Gamma_m=\Gamma(P_m, Q_m, D),$ тогда функция
$\rho(x)= \left\{
\begin{array}{rr}
\frac{1}{l}, & x\in D\\
0,  &  x\notin  D
\end{array}
\right. $
является допустимой для семейства $\Gamma_m,$ поскольку для
произвольной (локально спрямляемой) кривой $\gamma\in \Gamma_m$
выполнено $\int\limits_{\gamma}\rho(x)|dx|\geqslant
\frac{l(\gamma)}{l}\geqslant 1$ (где $l(\gamma)$ обозначает длину
кривой $\gamma$). Поскольку по условию отображения $f_m$
удовлетворяют (\ref{eq1C}), получаем:
\begin{equation}\label{eq14***}
M_{\alpha^{\,\prime}}(f_m(\Gamma_m))\leqslant
\frac{1}{l^{\alpha}}\,\,\int\limits_{D} Q(x)\,d\mu(x):=c<\infty\,,
\end{equation}
т.к. $Q\in L^1(D).$
С другой стороны, ${\rm diam}\,f_m(P_m)\geqslant d^{\,\prime}(y_m,
f_m(z^1_m)) \geqslant (1/2)\cdot d^{\,\prime}(y_0, p_1)>0$ и ${\rm
diam}\,f_m(Q_m)\geqslant d^{\,\prime}(y_0, f_m(z^2_m)) \geqslant
(1/2)\cdot d^{\,\prime}(y_0, p_2)>0,$ кроме того,
${\rm dist}\,(f_m(P_m), f_m(Q_m))\leqslant d^{\,\prime}(y_m,
y_0)\rightarrow 0,$ $m\rightarrow \infty\,.$ Тогда ввиду
пункта~\textbf{1)} леммы~\ref{lem2}
$$M_{\alpha^{\,\prime}}(f_m(\Gamma_m))=M_{\alpha^{\,\prime}}(f_m(P_m), f_m(Q_m), D^{\,\prime})\rightarrow\infty\,,\quad m\rightarrow\infty\,,$$
что противоречит соотношению (\ref{eq14***}). Полученное
противоречие указывает на ошибочность предположения в
(\ref{eq13***}), что и завершает доказательство теоремы.~$\Box$
\end{proof}

Для числа $\delta>0,$ областей $D\subset X,$ $D^{\,\prime}\subset
X^{\,\prime},$ континуума $A\subset D$ и произвольной измеримой по
Лебегу функции $Q(x): X\rightarrow [1, \infty],$ $Q(x)\equiv 0$ при
$x\not\in D,$ обозначим через ${\frak H}_{\delta, A, Q }(D,
D^{\,\prime})$ семейство всех $Q$-гомеоморфизмов $f:D\rightarrow
D^{\,\prime},$ $f(D)=D^{\,\prime},$ таких что ${\rm
diam}\,f(A)\geqslant\delta.$ Справедливо следующее утверждение.

\medskip
\begin{theorem}\label{th4} {\sl В условиях теоремы~\ref{th4A}
каждый элемент $g$ семейства
$${\frak H}^{\,-1}_{\delta, A, Q }(D,
D^{\,\prime}):=\left\{g=f^{\,-1}:D^{\,\prime}\rightarrow D,\,\,\,
f\in{\frak H}_{\delta, A, Q }(D, D^{\,\prime}) \right\}$$
может быть продолжен по непрерывности до отображения
$\overline{g}=\overline{f^{\,-1}}:
\overline{D^{\,\prime}}\rightarrow\overline{D},$ причём
$g(\overline{D^{\,\prime}})=\overline{D}$ и семейство ${\frak
H}^{\,-1}_{\delta, A, Q }(\overline{D},
\overline{D^{\,\prime}}):=\{\overline{g}:\overline{D^{\,\prime}}\rightarrow
\overline{D},\,\overline{g}|_{D^{\,\prime}}=g,\, g\in{\frak
H}^{\,-1}_{\delta, A, Q }(D, D^{\,\prime})\}$ является равностепенно
непрерывным в $\overline{D^{\,\prime}}.$}
\end{theorem}

\medskip
\begin{proof}
Возможность непрерывного продолжения каждого ${\frak
H}^{\,-1}_{\delta, A, Q }(D, D^{\,\prime})$ на границу области
$D^{\,\prime}$ есть утверждение леммы~\ref{lem2}, пункт~{\textbf
3)}, а равностепенная непрерывность семейства в ${\frak
H}^{\,-1}_{\delta, A, Q }(\overline{D}, \overline{D^{\,\prime}})$ в
$D^{\,\prime}$ вытекает из теоремы~\ref{th4A}. Равенство
$g(\overline{D^{\,\prime}})=\overline{D}$ для $g\in{\frak
H}^{\,-1}_{\delta, A, Q }$ является утверждением пункта~\textbf{4}
леммы~\ref{lem2}. Осталось показать равностепенную непрерывность
семейства ${\frak H}^{\,-1}_{\delta, A, Q }(\overline{D},
\overline{D^{\,\prime}})$ на границе области $D^{\,\prime}.$

Предположим противное, а именно, допустим, что найдётся точка
$z_0\in
\partial D^{\,\prime},$ число $\varepsilon_0>0$ и последовательности
$z_m\in \overline{D^{\,\prime}},$ $z_m\rightarrow z_0$ при
$m\rightarrow\infty$ и $\overline{g}_m\in {\frak H}^{\,-1}_{\delta,
A, Q }(\overline{D}, \overline{D^{\,\prime}})$ такие, что
\begin{equation}\label{eq12}
d(\overline{g}_m(z_m),\overline{g}_m(z_0))\geqslant\varepsilon_0,\quad
m=1,2,\ldots .
\end{equation}
Положим $g_m:=\overline{g}_m|_{D^{\,\prime}}.$ Так как $g_m$ по
непрерывности продолжается на границу $D^{\,\prime},$ можно считать,
что $z_m\in D^{\,\prime}$ и, значит, $\overline{g}_m(z_m)=g_m(z_m).$
Кроме того, найдётся ещё одна последовательность $z^{\,\prime}_m\in
D^{\,\prime},$ $z^{\,\prime}_m\rightarrow z_0$ при
$m\rightarrow\infty,$ такая, что
$d(g_m(z^{\,\prime}_m),\overline{g}_m(z_0))\rightarrow 0$ при
$m\rightarrow\infty.$ Так как $\overline{D}$ -- компакт, мы можем
считать, что последовательности $g_m(z_m)$ и $g_m(z_0)$ являются
сходящимися при $m\rightarrow\infty.$ Пусть $g_m(z_m)\rightarrow
\overline{x_1}$ и $\overline{g}_m(z_0)\rightarrow \overline{x_2}$
при $m\rightarrow\infty.$ По непрерывности модуля из (\ref{eq12})
вытекает, что $d(\overline{x_1},
\overline{x_2})\geqslant\varepsilon_0,$ более того, так как
гомеоморфизмы сохраняют границу, $\overline{x_2}\in\partial D.$
Пусть $x_1$ и $x_2$ -- произвольные различные точки континуума $A,$
ни одна из которых не совпадает с $\overline{x_1}.$ Ввиду
условия~\textbf{A} можно соединить точки $x_1$ и $\overline{x_1}$
кривой $\gamma_1:[0, 1]\rightarrow \overline{D},$ а точки $x_2$ и
$\overline{x_2}$ -- кривой $\gamma_2:[0, 1]\rightarrow \overline{D}$
так, что $|\gamma_1|\cap |\gamma_2|=\varnothing,$ $\gamma_i(t)\in D$
при всех $t\in (0, 1),$ $i=1,2,$ $\gamma_1(0)=x_1,$
$\gamma_1(1)=\overline{x_1},$ $\gamma_2(0)=x_2$ и
$\gamma_2(1)=\overline{x_2}.$ Так как $D$ локально связна на своей
границе, найдутся непересекающиеся окрестности $U_1$ и $U_2$ точек
$\overline{x_1}$ и $\overline{x_2},$ соответственно, такие, что
$W_i:=D\cap U_i$ -- линейно связное множество. За счёт уменьшения
окрестностей $U_i,$ если это необходимо, мы можем считать, что
$\overline{U_1}\cap\overline{U_2}=\varnothing$ и
$\overline{U_1}\cap|\gamma_2|=\varnothing=\overline{U_2}\cap|\gamma_1|.$
Мы также можем считать, что $g_m(z_m)\in W_1$ и
$g_m(z^{\,\prime}_m)\in W_2$ при всех $m\in {\Bbb N}.$ Пусть $a_1$ и
$a_2$ -- произвольные точки, принадлежащие $|\gamma_1|\cap W_1$ и
$|\gamma_2|\cap W_2.$ Пусть $t_1, t_2$ таковы, что
$\gamma_1(t_1)=a_1$ и $\gamma_2(t_2)=a_2.$ Соединим точку $a_1$ с
точкой $g_m(z_m)$ кривой $\alpha_m:[t_1, 1]\rightarrow W_1$ такой,
что $\alpha_m(t_1)=a_1$ и $\alpha_m(1)=g_m(z_m).$ Аналогично,
соединим точку $a_2$ с точкой $g_m(z^{\,\prime}_m)$ кривой
$\beta_m:[t_2, 1]\rightarrow W_2$ такой, что $\beta_m(t_2)=a_2$ и
$\beta_m(1)=g_m(z^{\,\prime}_m)$ (см. рисунок~\ref{fig3}).
\begin{figure}[h]
\centerline{\includegraphics[scale=0.6]{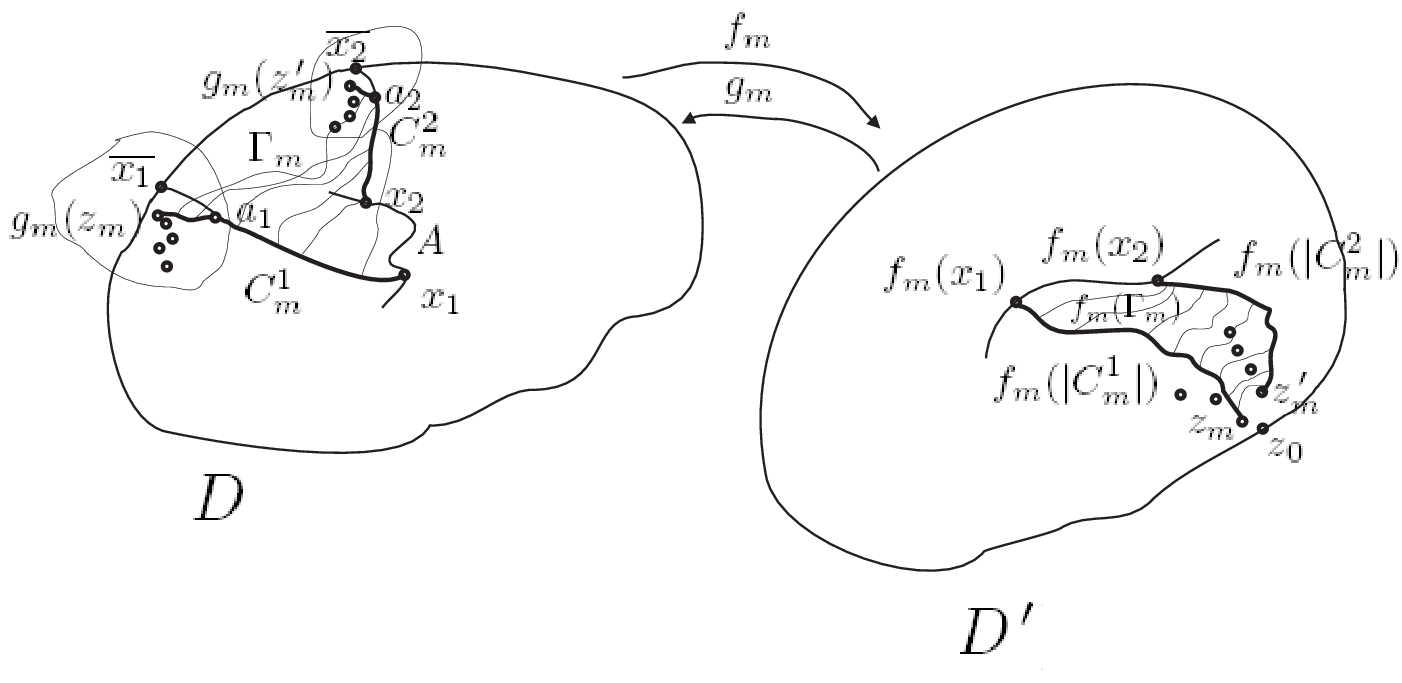}} \caption{К
доказательству теоремы~\ref{th4}}\label{fig3}
\end{figure}
\medskip
Положим теперь
$C^1_m(t)=\quad\left\{
\begin{array}{rr}
\gamma_1(t), & t\in [0, t_1],\\
\alpha_m(t), & t\in [t_1, 1]\end{array} \right.\,,$
$C^2_m(t)=\quad\left\{
\begin{array}{rr}
\gamma_2(t), & t\in [0, t_2],\\
\beta_m(t), & t\in [t_2, 1]\end{array} \right.\,.$
Пусть, как обычно, $|C^1_m|$ и $|C^2_m|$ -- носители кривых $C^1_m$
и $C^2_m,$ соответственно. Заметим, что по построению $|C^1_m|$ и
$|C^2_m|$ -- два непересекающихся континуума в $D,$ причём
существует $l_0>0$ такое, что ${\rm dist}\,(|C^1_m|, |C^2_m|)>l_0>0$
при всех $m=1,2,\ldots .$ Пусть теперь $\Gamma_m$ -- семейство
кривых, соединяющих $|C^1_m|$ и $|C^2_m|$ в $D.$ Тогда функция
$\rho(x)= \left\{
\begin{array}{rr}
\frac{1}{l_0}, & x\in D\\
0,  &  x\notin  D
\end{array}
\right.$
является допустимой для семейства $\Gamma_m$ и, поскольку по условию
отображения $f_m,$ $f_m=g_m^{\,-1},$ удовлетворяют (\ref{eq1C}),
получаем:
%
$$M_{\alpha^{\,\prime}}(f_m(\Gamma_m))\leqslant
\frac{1}{l_0^{\alpha}}\,\,\int\limits_{D} Q(x)\,d\mu(x):=c=c(l_0,
Q)<\infty\,,$$
т.к. $Q\in L^1(D).$
С другой стороны, ввиду пункта~\textbf{5} леммы~\ref{lem2} найдётся
число $\delta_1>0$ такое, что ${\rm dist}\,(f_{m}(A), \partial
D^{\,\prime})>\delta_1>0,$ $m=1,2,\ldots \,.$ Отсюда получим, что
найдётся некоторый номер $m_1>m_0,$ $m_0\in {\Bbb N},$ такой что при
всех $m\geqslant m_1$
$${\rm diam}\,f_m(|C^1_m|)\geqslant d^{\,\prime}(z_m, f_m(x_1)) \geqslant
(1/2)\cdot {\rm dist}\,(f_m(A), \partial
D^{\,\prime})>\delta_1/2\,,$$
\begin{equation}\label{eq14}
{\rm diam}\,f_m(|C^2_m|)\geqslant d^{\,\prime}(z^{\,\prime}_m,
f_m(x_2)) \geqslant (1/2)\cdot {\rm dist}\,(f_m(A), \partial
D^{\,\prime})>\delta_1/2\,.
\end{equation}
Кроме того, ${\rm dist}\,(f_m(|C^1_m|, f_m(|C^2_m|)\leqslant
d^{\,\prime}(z_m, z^{\,\prime}_m)\rightarrow 0$ по выбору $z_m$ и
$z^{\,\prime}_m.$ Тогда в силу пункта~\textbf{1} леммы~\ref{lem2}
$M_{\alpha^{\,\prime}}(f_m(\Gamma_m))=M_{\alpha^{\,\prime}}(\Gamma(f_m(|C^1_m|),
f_m(|C^2_m|), D^{\,\prime}))\rightarrow\infty$ при
$m\rightarrow\infty,$ но это противоречит соотношению~(\ref{eq14}).
Полученное противоречие указывает на неверность предположения
в~(\ref{eq12}), что и доказывает теорему.~$\Box$
\end{proof}

\medskip Для числа $\delta>0,$ областей $D, D^{\,\prime}\subset {\Bbb R}^n,$ $n\geqslant 2,$
континуума $A\subset D$ и произвольной измеримой по Лебегу функции
$Q(x): {\Bbb R}^n\rightarrow [1, \infty],$ $Q(x)\equiv 0$ при
$x\not\in D,$ обозначим через ${\frak S}_{\delta, A, Q }(D,
D^{\,\prime})$ семейство всех $Q$-гомеоморфизмов $f:D\rightarrow
D^{\,\prime},$ $f(D)=D^{\,\prime},$ таких что ${\rm
diam}\,f(A)\geqslant\delta.$ В качестве простого следствия из
теоремы \ref{th4} получаем следующее утверждение для пространства
${\Bbb R}^n$ (см. также~\cite{Sev$_3$}).

\medskip
\begin{corollary}\label{cor3} {\sl Предположим, что
область $D$ локально связна во всех граничных точках, $\overline{D}$
и $\overline{D^{\,\prime}}$ являются компактами, $D^{\,\prime}$
является $QED$-об\-лас\-тью и $Q\in L^1(D).$
Тогда каждый элемент $g$ семейства ${\frak S}^{\,-1}_{\delta, A, Q
}(D, D^{\,\prime}),$ состоящего из всех обратных гомеоморфизмов
$\left\{g=f^{-1}:D^{\,\prime}\rightarrow D| f\in {\frak S}_{\delta,
A, Q }(D, D^{\,\prime})\right\}$ может быть продолжен по
непрерывности до отображения $\overline{g}=\overline{f^{\,-1}}:
\overline{D^{\,\prime}}\rightarrow\overline{D},$ причём семейство
${\frak S}^{\,-1}_{\delta, A, Q }(\overline{D},
\overline{D^{\,\prime}}),$ состоящее из всех продолженных, таким
образом, отображений является равностепенно непрерывным в
$\overline{D^{\,\prime}}.$}
\end{corollary}

\medskip
Заметим, что при $Q(x)=const$ утверждение следствия~\ref{cor3} есть
ослабленный вариант теоремы Някки-Палка о равностепенной
непрерывности семейства квазиконформных отображений в замыкании
области (см.~\cite[теорема~3.1]{NP}).

\medskip
{\bf 5. О выполнении условия \textbf{A} в евклидовом пространстве.}
Для открытого, замкнутого, либо полуоткрытого интервала $I\subset
{\Bbb R}$ и кривой $\gamma: I\rightarrow {\Bbb R}^n,$ как обычно,
полагаем $|\gamma|=\{x\in {\Bbb R}^n: \exists\, t\in [a, b]:
\gamma(t)=x\}.$ Как известно, кривая $\gamma:I\rightarrow {\Bbb
R}^n$ называется {\it жордановой дугой}, если $\gamma$ --
гомеоморфизм на $I.$ Справедливо следующее утверждение.

\medskip
\begin{proposition}\label{pr1}
Пусть $D$ -- область в ${\Bbb R}^n,$ $n\geqslant 2,$ локально
связная на своей границе. Тогда любые две пары точек $a\in D, b\in
\overline{D},$ и $c\in D, d\in \overline{D}$ можно соединить
непересекающимися между собой кривыми $\gamma_1:[0, 1]\rightarrow
\overline{D}$ и $\gamma_2:[0, 1]\rightarrow \overline{D},$  такими,
что $\gamma_i(t)\in D$ при всех $t\in (0, 1),$ $i=1,2,$
$\gamma_1(0)=a,$ $\gamma_1(1)=b,$ $\gamma_2(0)=c,$ $\gamma_2(1)=d.$
\end{proposition}

\medskip
\begin{proof}
Заметим, что точки области, локально связной на границе, являются
достижимыми изнутри области посредством кривых (см.
\cite[предложение~13.2]{MRSY}). В таком случае, если $n\geqslant 3,$
соединим точки $a$ и $b$ произвольной жордановой дугой $\gamma_1$ в
области $D,$ не проходящей через точки $c$ и $d$ (что возможно по
ввиду локальной связности $D$ на границе и переходом от кривой к
ломаной, если это необходимо). Тогда $\gamma_1$ не разбивает область
$D$ как множество топологической размерности 1 (см.
\cite[следствие~1.5.IV]{HW}), что и обеспечивает существование
искомой кривой $\gamma_2.$ Таким образом, в случае $n\geqslant 3$
утверждение леммы~\ref{lem1} установлено.

Пусть теперь $n=2,$ тогда снова точки $c$ и $d$ не разбивают область
$D$ (\cite[следствие~1.5.IV]{HW}). В таком случае, также можно
соединить точки $a$ и $b$ жордановой дугой $\gamma_1$ в $D,$ не
проходящей через точки $c$ и $d.$ Ввиду теоремы Антуана (см.
\cite[теорема~4.3, \S\,4]{Keld}) область $D$ можно отобразить на
некоторую область $D^{\,*}$ посредством плоского гомеоморфизма
$\varphi:{\Bbb R}^2\rightarrow {\Bbb R}^2$ так, что
$\varphi(\gamma_1)=J$ и $J$ -- отрезок в $D^{\,*}.$ Заметим также,
что точки границы области $D^{\,*}$ являются достижимыми изнутри
$D^{\,*}$ посредством кривых. Таким образом, мы можем соединить
точки $\varphi(c)$ и $\varphi(d)$ в $D^{\,*}$ посредством жордановой
кривой $\alpha_2:[0, 1]\rightarrow \overline{D^{\,*}},$ которая
целиком лежит в $D^{\,*},$ кроме, может быть, своей концевой точки
$\alpha_2(1)=\varphi(d).$

Осталось показать, что кривую $\alpha_2$ можно выбрать так, что она
не будет пересекать отрезок $J.$ В самом деле, пусть $\alpha_2$
пересекает $J,$ и пусть $t_1$ и $t_2$ -- соответственно, наибольшее
и наименьшее значение $t\in [0, 1],$ для которых $\alpha_2(t)\in
|J|.$ Пусть также
$$J=J(s)=\varphi(a)+ (\varphi(b)-\varphi(a))s, \quad s\in [0, 1]$$
-- параметризация отрезка $J.$ Пусть $\widetilde{s_1}$ и
$\widetilde{s_2}\in (0, 1)$ таковы, что
$J(\widetilde{s_1})=\alpha_2(t_1)$ и
$J(\widetilde{s_2})=\alpha_2(t_2).$ Положим
$s_2:=\max\{\widetilde{s_1}, \widetilde{s_2}\}.$ Пусть
$e_1=\varphi(b)-\varphi(a)$ и $e_2$ -- единичный вектор,
ортогональный $e_1,$ тогда множество
$$P_{\varepsilon}=\{x=\varphi(a)+
x_1e_1+x_2e_2,\quad x_1\in (-\varepsilon, s_2+\varepsilon), \quad
x_2\in (-\varepsilon, \varepsilon)\}\,,\quad \varepsilon>0\,,$$
представляет собой прямоугольник, содержащий $|J_1|,$ где $J_1$ --
сужение $J$ на отрезок $[0, s_2]$ (см. рисунок~\ref{fig4}).
\begin{figure}[h]
\centerline{\includegraphics[scale=0.5]{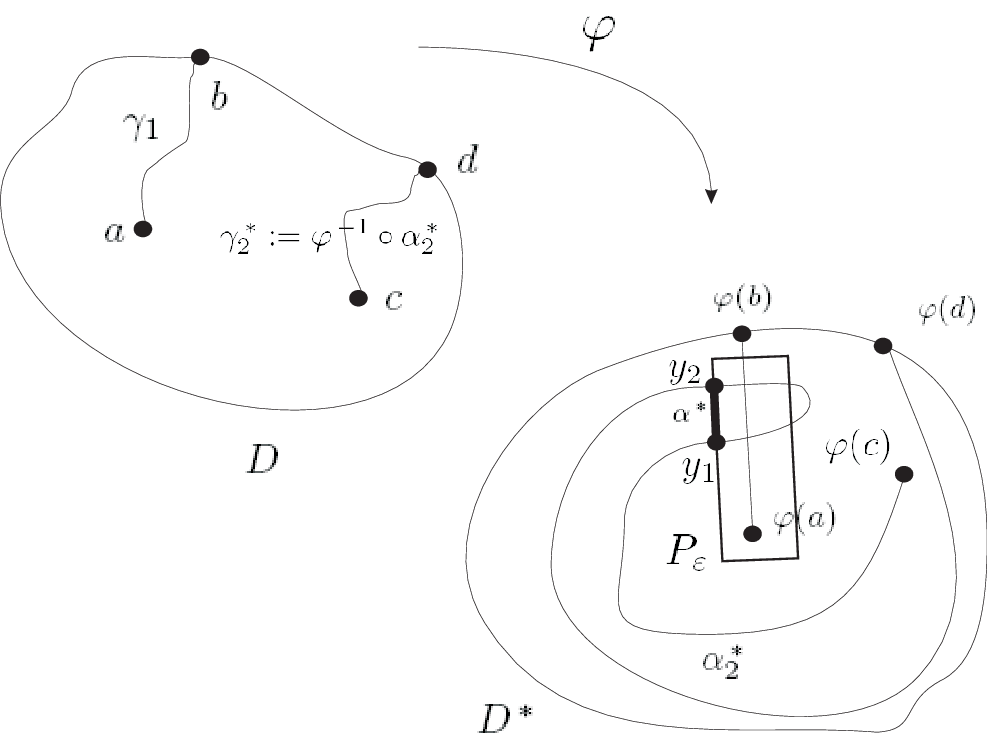}}
\caption{Возможность соединения двух пар точек кривыми в области
}\label{fig4}
\end{figure}
Выберем $\varepsilon>0$ так, что $\varphi(c)\not\in
P_{\varepsilon},$ ${\rm dist}\,(P_{\varepsilon}, \partial
D^{\,*})>\varepsilon.$ Ввиду \cite[теорема~1.I, гл.~5, \S\, 46]{Ku})
кривая $\alpha_2$ пересекает $\partial P_{\varepsilon}$ при
некотором $T_1<t_1$ и при некотором $T_2>t_2.$ Пусть
$\alpha_2(T_1)=y_1$ и $\alpha_2(T_2)=y_2.$ Так как $\partial
P_{\varepsilon}$ -- связное множество, можно соединить точки $y_1$ и
$y_2$ кривой $\alpha^{\,*}(t):[T_1, T_2]\rightarrow
\partial P_{\varepsilon}.$ Окончательно, положим
$$\alpha_2^{\,*}(t)\quad =\quad\left\{
\begin{array}{rr}
\alpha_2(t), & t\in [0, 1]\setminus [T_1, T_2],\\
\alpha^{\,*}(t), & t\in [T_1, T_2]\end{array} \right.$$
и $\gamma^{\,*}_2:=\varphi^{\,-1}\circ \alpha_2^{\,*}.$ Тогда
$\gamma_1$ соединяет $a$ и $b$ в $D,$ а $\gamma_2^{\,*}$ соединяет
$c$ и $d$ в $D,$ при этом, $\gamma_1$ и $\gamma_2^{\,*}$ не
пересекаются, что и следовало установить.
\end{proof}

\medskip
{\bf Пример~1.} Как известно, дробно-линейные автоморфизмы
единичного круга ${\Bbb D}\subset{\Bbb C}$ на себя задаются формулой
$f(z)=e^{i\theta}\frac{z-a}{1-\overline{a}z},$ $z\in {\Bbb D},$
$a\in{\Bbb C},$ $|a|<1,$ $\theta\in [0, 2\pi).$ Указанные
отображения $f$ являются 1-гомеоморфизмами; все условия
следствия~\ref{cor3} выполняются, кроме условия ${\rm
diam}\,f(A)\geqslant\delta,$ выполнение которого зависит от
конкретного вида этих отображений.

Если, например, $\theta=0$ и $a=1/n,$ $n=1,2,\ldots,$ то
$f_n(z)=\frac{z-1/n}{1-z/n}=\frac{nz-1}{n-z}.$ Положим $A=[0, 1/2],$
тогда $f_n(0)=-1/n\rightarrow 0$ и
$f_n(1/2)=\frac{n-2}{2n-1}\rightarrow 1/2,$ $n\rightarrow\infty.$
Отсюда видно, что последовательность $f_n$ удовлетворяет условию
${\rm diam}\,f_n(A)\geqslant\delta,$ например, при $\delta=1/4.$
Путём непосредственных вычислений убеждаемся в том, что
$f_n^{\,-1}(z)=\frac{z+1/n}{1+z/n}$ и, значит, $f_n^{\,-1}$
равномерно сходится к $f^{\,-1}(z)\equiv z.$ Таким образом,
последовательность $f_n^{\,-1}(z)$ равностепенно непрерывна в
$\overline{{\Bbb D}}$ (что, впрочем, независимо от прямых вычислений
вытекает из следствия~\ref{cor3}).

\medskip
Если же взять
$f^{-1}_n(z)=\frac{z-(n-1)/n}{1-z(n-1)/n}=\frac{nz-n+1}{n-nz+1},$
то, как легко видеть, такая последовательность локально равномерно
({\it но не равномерно !}) сходится к $-1$ в ${\Bbb D};$ в то же
время, $f^{\,-1}_n(1)=1.$ Следовательно, $f^{\,-1}_n$ не является
равностепенно непрерывной в точке 1; в этом случае
$f_n(z)=\frac{z+(n-1)/n}{1+z(n-1)/n}$ и условие ${\rm
diam}\,f_n(A)\geqslant\delta$ ни при каком $\delta>0,$ не зависящем
от $n,$ не может быть выполнено ввиду следствия~\ref{cor3} (что,
впрочем, непосредственно видно, поскольку $f_n(z)\rightarrow 1$ при
$n\rightarrow \infty$ локально равномерно, но не равномерно, в
${\Bbb D}$).

\medskip
{\bf Пример~2.} Пусть $p\geqslant 1$ настолько велико, что число
$n/p(n-1)$ меньше 1, и пусть, кроме того $\alpha\in (0, n/p(n-1))$
-- произвольное число. Рассмотрим последовательность отображений
$f_m: {\Bbb B}^n\rightarrow B(0, 2)$ шара ${\Bbb B}^n$ на шар $B(0,
2)$ следующим образом:
$$f_m(x)\,=\,\left
\{\begin{array}{rr} \frac{1+|x|^{\alpha}}{|x|}\cdot x\,, & 1/m\leqslant|x|\leqslant 1, \\
\frac{1+(1/m)^{\alpha}}{(1/m)}\cdot x\,, & 0<|x|< 1/m \ .
\end{array}\right.
$$
Заметим, что $f_m$ являются $Q$-гомеоморфизмами в ${\Bbb B}^n$ при
$Q=\left(\frac{1+|x|^{\,\alpha}}{\alpha
|x|^{\,\alpha}}\right)^{n-1}\in L^1({\Bbb B}^n)$
(см.~\cite[доказательство теоремы~7.1]{Sev$_3$}), что $B(0, 2)$
является $QED$-областью (см.~\cite[лемма~4.3]{Vu}) и что ${\Bbb
R}^n$ является пространством Лёвнера (см.~\cite[теорема~8.2]{He}).
По построению отображения $f_m$ фиксируют бесконечное число точек
единичного шара при всех $m\geqslant 2.$ Равностепенная
непрерывность семейства отображений $g_m:=f_m^{\,-1}$ в $B(0, 2),$
$$g_m(y):=f^{-1}_m(y)\,=\,\left
\{\begin{array}{rr} \frac{y}{|y|}(|y|-1)^{1/\alpha}\,, & 1+1/m^{\alpha}\leqslant|y|< 2, \\
\frac{(1/m)}{1+(1/m)^{\alpha}}\cdot y\,, & 0<|y|< 1+1/m^{\alpha} \ ,
\end{array}\right.
$$
вытекает прямо из теоремы~\ref{th4}, однако, может быть установлена
и непосредственно. Из этой же теоремы вытекает равностепенная
непрерывность продолженного по непрерывности семейства
$\{g_m\}_{m=1}^{\infty}$ на $\overline{B(0, 2)}.$

\medskip
Обращает на себя внимание то обстоятельство, что хотя семейство
отображений $\frak G=\{g_m\}_{m=1}^{\infty}$ равностепенно
непрерывно в $B(0, 2)$, таковым не является <<обратное>> к ним
семейство $\frak F=\{f_m\}_{m=1}^{\infty}$ (в самом деле,
$|f_m(x_m)-f(0)|=1+1/m\not\rightarrow 0$ при $m\rightarrow\infty,$
где $|x_m|=1/m$).

Теперь исследуем последовательность $g_m$ в контексте утверждения
теоремы \ref{th3}. Нетрудно проверить, что последовательность $g_m$
равномерно сходится в $B(0, 2)$ к отображению
$$g(y):=\,\left
\{\begin{array}{rr} \frac{y}{|y|}(|y|-1)^{1/\alpha}\,, & 1<|y|< 2, \\
0\,, & 0<|y|< 1+1/m^{\alpha} \ ,
\end{array}\right.
$$
которое не является ни постоянным, ни нульмерным. В силу теоремы~7.1
в \cite{Sev$_3$}, заключаем, что $g_m$ удовлетворяет условию
(\ref{eq2A}) (и, тем более, условию (\ref{eq2*A}) при $p=q=n$) при
$Q=\left(\frac{1+|x|^{\,\alpha}}{\alpha
|x|^{\,\alpha}}\right)^{n-1}\in L^p,$ где $p$ -- число, выбранное в
начале рассмотрения данного примера. Последнее указывает на то, что
условия на функцию $Q,$ присутствующие в теореме \ref{th3}, являются
точными в следующем смысле: требование $Q\in FMO$ нельзя заменить на
условие $Q\in L^p$ ни для какого (сколь угодно большого $p>1$).

КОНТАКТНАЯ ИНФОРМАЦИЯ

\medskip
\noindent{{\bf Евгений Александрович Севостьянов} \\
\noindent{{\bf Сергей Александрович Скворцов} \\
Житомирский государственный университет им.\ И.~Франко\\
кафедра математического анализа, ул. Большая Бердичевская, 40 \\
г.~Житомир, Украина, 10 008 \\ тел. +38 066 959 50 34 (моб.),
e-mail: esevostyanov2009@gmail.com}

\end{document}